\documentclass{amsart}
\usepackage{amssymb,amsmath, euscript}

\newcounter{sec}

\def\SS{\smallskip}

\newcounter{punct}[sec]

\def\punct{\refstepcounter{punct}{\arabic{sec}.\arabic{punct}.  }}

\def\COUNTERS{\addtocounter{sec}{1}
              \setcounter{punct}{0}
          \setcounter{equation}{0}
          \setcounter{theorem}{0}
                  }

\newtheorem{theorem}{Theorem}[sec]
\newtheorem{proposition}[theorem]{Proposition}

\begin{document}

\def\ov{\overline}

\def\wt{\widetilde}

\newcommand{\rk}{\mathop {\mathrm {rk}}\nolimits}
\newcommand{\Aut}{\mathop {\mathrm {Aut}}\nolimits}
\newcommand{\Out}{\mathop {\mathrm {Out}}\nolimits}
\renewcommand{\Re}{\mathop {\mathrm {Re}}\nolimits}

\def\Br{\mathrm {Br}}

\def\SL{\mathrm {SL}}
\def\SU{\mathrm {SU}}
\def\GL{\mathrm  {GL}}
\def\U{\mathrm  U}
\def\OO{\mathrm  O}
\def\Sp{\mathrm  {Sp}}
\def\SO{\mathrm  {SO}}
\def\SOS{\mathrm {SO}^*}
\def\Diff{\mathrm{Diff}}
\def\Vect{\mathfrak{Vect}}

\def\PGL{\mathrm  {PGL}}
\def\PU{\mathrm {PU}}

\def\PSL{\mathrm  {PSL}}

\def\Symp{\mathrm{Symp}}

\def\cA{\mathcal A}
\def\cB{\mathcal B}
\def\cC{\mathcal C}
\def\cD{\mathcal D}
\def\cE{\mathcal E}
\def\cF{\mathcal F}
\def\cG{\mathcal G}
\def\cH{\mathcal H}
\def\cJ{\mathcal J}
\def\cI{\mathcal I}
\def\cK{\mathcal K}
\def\cL{\mathcal L}
\def\cM{\mathcal M}
\def\cN{\mathcal N}
\def\cO{\mathcal O}
\def\cP{\mathcal P}
\def\cQ{\mathcal Q}
\def\cR{\mathcal R}
\def\cS{\mathcal S}
\def\cT{\mathcal T}
\def\cU{\mathcal U}
\def\cV{\mathcal V}
\def\cW{\mathcal W}
\def\cX{\mathcal X}
\def\cY{\mathcal Y}
\def\cZ{\mathcal Z}



\def\frA{\mathfrak A}
\def\frB{\mathfrak B}
\def\frC{\mathfrak C}
\def\frD{\mathfrak D}
\def\frE{\mathfrak E}
\def\frF{\mathfrak F}
\def\frG{\mathfrak G}
\def\frH{\mathfrak H}
\def\frJ{\mathfrak J}
\def\frK{\mathfrak K}
\def\frL{\mathfrak L}
\def\frM{\mathfrak M}
\def\frN{\mathfrak N}
\def\frO{\mathfrak O}
\def\frP{\mathfrak P}
\def\frQ{\mathfrak Q}
\def\frR{\mathfrak R}
\def\frS{\mathfrak S}
\def\frT{\mathfrak T}
\def\frU{\mathfrak U}
\def\frV{\mathfrak V}
\def\frW{\mathfrak W}
\def\frX{\mathfrak X}
\def\frY{\mathfrak Y}
\def\frZ{\mathfrak Z}

\def\fra{\mathfrak a}
\def\frb{\mathfrak b}
\def\frc{\mathfrak c}
\def\frd{\mathfrak d}
\def\fre{\mathfrak e}
\def\frf{\mathfrak f}
\def\frg{\mathfrak g}
\def\frh{\mathfrak h}
\def\fri{\mathfrak i}
\def\frj{\mathfrak j}
\def\frk{\mathfrak k}
\def\frl{\mathfrak l}
\def\frm{\mathfrak m}
\def\frn{\mathfrak n}
\def\fro{\mathfrak o}
\def\frp{\mathfrak p}
\def\frq{\mathfrak q}
\def\frr{\mathfrak r}
\def\frs{\mathfrak s}
\def\frt{\mathfrak t}
\def\fru{\mathfrak u}
\def\frv{\mathfrak v}
\def\frw{\mathfrak w}
\def\frx{\mathfrak x}
\def\fry{\mathfrak y}
\def\frz{\mathfrak z}

\def\frsp{\mathfrak{sp}}


\def\bfa{\mathbf a}
\def\bfb{\mathbf b}
\def\bfc{\mathbf c}
\def\bfd{\mathbf d}
\def\bfe{\mathbf e}
\def\bff{\mathbf f}
\def\bfg{\mathbf g}
\def\bfh{\mathbf h}
\def\bfi{\mathbf i}
\def\bfj{\mathbf j}
\def\bfk{\mathbf k}
\def\bfl{\mathbf l}
\def\bfm{\mathbf m}
\def\bfn{\mathbf n}
\def\bfo{\mathbf o}
\def\bfp{\mathbf p}
\def\bfq{\mathbf q}
\def\bfr{\mathbf r}
\def\bfs{\mathbf s}
\def\bft{\mathbf t}
\def\bfu{\mathbf u}
\def\bfv{\mathbf v}
\def\bfw{\mathbf w}
\def\bfx{\mathbf x}
\def\bfy{\mathbf y}
\def\bfz{\mathbf z}

\def\bfA{\mathbf A}
\def\bfB{\mathbf B}
\def\bfC{\mathbf C}
\def\bfD{\mathbf D}
\def\bfE{\mathbf E}
\def\bfF{\mathbf F}
\def\bfG{\mathbf G}
\def\bfH{\mathbf H}
\def\bfI{\mathbf I}
\def\bfJ{\mathbf J}
\def\bfK{\mathbf K}
\def\bfL{\mathbf L}
\def\bfM{\mathbf M}
\def\bfN{\mathbf N}
\def\bfO{\mathbf O}
\def\bfP{\mathbf P}
\def\bfQ{\mathbf Q}
\def\bfR{\mathbf R}
\def\bfS{\mathbf S}
\def\bfT{\mathbf T}
\def\bfU{\mathbf U}
\def\bfV{\mathbf V}
\def\bfW{\mathbf W}
\def\bfX{\mathbf X}
\def\bfY{\mathbf Y}
\def\bfZ{\mathbf Z}

\def\bfw{\mathbf w}

\def\R {{\mathbb R }}
 \def\C {{\mathbb C }}
  \def\Z{{\mathbb Z}}
  \def\H{{\mathbb H}}
\def\K{{\mathbb K}}
\def\N{{\mathbb N}}
\def\Q{{\mathbb Q}}
\def\A{{\mathbb A}}

\def\T{\mathbb T}
\def\P{\mathbb P}

\def\bbA{\mathbb A}
\def\bbB{\mathbb B}
\def\bbD{\mathbb D}
\def\bbE{\mathbb E}
\def\bbF{\mathbb F}
\def\bbG{\mathbb G}
\def\bbI{\mathbb I}
\def\bbJ{\mathbb J}
\def\bbL{\mathbb L}
\def\bbM{\mathbb M}
\def\bbN{\mathbb N}
\def\bbO{\mathbb O}
\def\bbP{\mathbb P}
\def\bbQ{\mathbb Q}
\def\bbS{\mathbb S}
\def\bbT{\mathbb T}
\def\bbU{\mathbb U}
\def\bbV{\mathbb V}
\def\bbW{\mathbb W}
\def\bbX{\mathbb X}
\def\bbY{\mathbb Y}

\def\kappa{\varkappa}
\def\epsilon{\varepsilon}
\def\phi{\varphi}
\def\le{\leqslant}
\def\ge{\geqslant}

         \begin{center}
\bf\Large   Spectral data for pairs of matrices
of order 3 and action of the group
 $\GL(2,\Z)$

\bigskip
\sc\large{Yury Neretin}%
\footnote{Supported by grants   FWF, project P19064,
 NWO.047.017.015, JSPS-RFBR-07.01.91209.}

\end{center}

{\small
The group $\GL(2,\Z)$ acts in a natural way on the
set of pairs of $n\times n$-matrices determined up to
a simultaneous conjugation. For $n=3$ we write
explicit formulas for action of generators of
$\GL(2,\Z)$ in the terms of spectral data of matrices,
i.e., spectral curves and line bundles.
}

\section{Introduction}

\COUNTERS

{\bf\punct Spectral data for pairs of matrices.}
Let  $A$, $B$ be a pair of nondegenerate matrices of order
  $n\times n$ defined up to to a simultaneous conjugation,
$$
 (A, B)\sim(gAg^{-1},gBg^{-1})
.$$
Denote the quotient space by $\frM(n)$.

The equation
$$
  \det(\lambda+\mu A+\nu B)=0
$$
determines a curve
 $\frC$ of degree $n$ on the projective plane $\C\P^2$.
Obviously $\frC$ depends on
class of equivalence and not on the matrices
$A$, $B$
themselves. For matrices in a general position
this curve is irreducible and nondegenerate.

It is impossible to reconstruct  a point $\in\frM(n)$ from  the
curve $\frC$ in a unique way. Indeed, $\dim\frM(n)=n^2+1$, on the
other hand the dimension of the space of curves of degree $n$ on
$\C\P^2$ is $C_n^2-1=\frac 12 (n+1)(n+2)-1$.
 For  $n>2$ the number of coordinates is not sufficient.

Next, consider the subspace
$$
\ker (\lambda+\mu A+\nu B)
\subset \C^n
.
$$
For matrices  $A$, $B$ of general position
these kernels are one-dimensional
at all points of the curve  $\frC$.
Thus, we get a linear bundle on  $\frC$,
say
$\cL$. By $\cL^*$ we denote the dual linear bundle.

A pair $(A,B)\in\frM(n)$ in general position can be reconstructed
from the curve
 $\subset\C\P^2$
and the linear bundle on the curve.
 The following statement was discovered many times, see
e.g., \cite{Zlo}) \footnote{See also an old theorem of Tyurin
 \cite{Tyu}
on triples of quadratic forms, similar constructions are
widely explored in
theory of integrable systems, see, e.g., \cite{RS}.}.

\begin{theorem}
 The map $(A,B)\mapsto (\frC,\cL^*)$ is a bijection
(up to a set of zero measure) of the space  $\frM(n)$ to the space
of curves of degree  $n$ on $\C\P^2$ equipped with linear bundles of
degree $n(n-1)/2$.
\end{theorem}


\smallskip

{\bf\punct The group of outer automorphisms
of a free group.}
Denote by $F_2$ the free group with
2 generators
$c_1$, $c_2$. Denote by $\Aut(F_2)$ the group
of automorphisms of $F_2$, by
$$
\Out(F_2)=\Aut(F_2)/F_2
$$
the group of outer automorphisms of $F_2$.

To define an element of  $\Aut(F_2)$,
we must indicate images
  of generators $c_1$, $c_2$,
\begin{align}
c_1\mapsto c_1^{\alpha_1} c_2^{\beta_1} c_1^{\alpha_2} c_2^{\beta_2}\dots
\label{eq:obrobr-1}
\\
c_2\mapsto c_1^{\gamma_1} c_2^{\delta_1} c_1^{\gamma_2} c_2^{\delta_2}\dots
\label{eq:obrobr-2}
\end{align}

 Elements  $\in F_2$ in the right hand side
are not arbitrary.
By the well-known Nielsen theorem
  (see, e.g., \cite{LS}),
 the group $\Aut(F_2)$
is generated by the following transformations:

\SS

a) $(c_1, c_2)\mapsto (c_2,c_1)$;

\SS

b) $(c_1,c_2)\mapsto (c_1^{-1}, c_2)$;

\SS

c) $(c_1,c_2)\mapsto(c_1, c_1 c_2)$.

Consider the commutant  $[F_2,F_2]$ of the free group,
i.e., the subgroup generated by all
products of the type
$xyx^{-1} y^{-1}$. Clearly, this subgroup is
invariant with respect to  $\Aut(F_2)$.
The quotient group $F_2/[F_2,F_2]$
is isomorphic to
 $\Z\oplus\Z$. Therefore $\Aut(F_2)$
acts by automorphisms on  $\Z\oplus\Z$. In other words, we get a
homomorphism
$$
\sigma:\Aut(F_2)\to \GL(2,\Z)
,$$
where $\GL(2,\Z)$ is the group of $2\times 2$ matrices
 $g$ with integer elements such that
 $g^{-1}$ also has integer elements%
\footnote{Automatically, $\det g=\pm1$.}.
The images of the transformations  a), b), c)
are given by matrices
$$
  \begin{pmatrix}0&1\\1&0 \end{pmatrix},\qquad
  \begin{pmatrix}-1&0\\0&1 \end{pmatrix},\qquad
  \begin{pmatrix}1&1\\0&1 \end{pmatrix},
$$
the transformation
 (\ref{eq:obrobr-1})--(\ref{eq:obrobr-2}) corresponds to the matrix
$$
  \begin{pmatrix} \sum\alpha_j& \sum\beta_j
\\
\sum\gamma_j&\sum\delta_j \end{pmatrix}
.
$$
Obviously inner automorphisms of
$F_2$
are contained in the kernel of  $\sigma$,
i.e., we get a homomorphism
$$
\sigma:\Out(F_2)\to \GL(2,\Z)
.$$
It turns out to be that (see, e.g., \cite{LS})
that $\sigma$ is an isomorphism.


\smallskip

{\bf \punct  Action of the group $\GL(2,\Z)$ on $\frM(n)$.}
The group $\Out(F_2)$  acts on  $\frM(n)$
in the obvious way,
the transformations (\ref{eq:obrobr-1})--(\ref{eq:obrobr-2})
correspond to
\begin{align*}
A\mapsto A_1^{\alpha_1} B_2^{\beta_1} A_1^{\alpha_2}
 B_2^{\beta_2}\dots
 \,,
\\
B\mapsto A_1^{\gamma_1} B_2^{\delta_1} A_1^{\gamma_2} B_2^{\delta_2}\dots
\,.\end{align*}

Problem:  {\it
is it possible to write this action in the terms of spectral data?}

For finite collections of
 $2\times 2$ matrices this problem was investigated
in different (nonequivalent) versions,
see, e.g.,  \cite{Fock}, \cite{Ner}.
In this note, we solve it for pairs
of $3\times 3$  matrices.
Unlikely, our solution is final
(even for triples of matrices),
in any case explicit formulas are written.


\smallskip

{\bf\punct Preliminary description of construction.}
{\it In that follows, we discuss only  $3\times 3$
matrices of general position.}

 The equation
$\det(\lambda+\mu A+\nu B)=0$ determines
a cubic curve  $\frC$ on the projective plane  $\C\P^2$.
The curve is determined by its equation.
Knowing the curve, we can assume that we know its intersections with
the line  $\nu=0$.

It is more difficult to introduce a convenient coordinate on the
space of linear bundles on $\frC$. Recall that a linear bundle is
determined by a class of equivalent divisors, see, e.g.,
 \cite{GH}. For a curve, a divisor
is a finite collection  $\{k_j x_j\}$ of points $x_j$
of the curve with integer
(positive or negative)
multiplicities $k_j$.
These collections arise as collections
of zeros and poles of meromorphic sections
of  line bundles.
By definition, {\it the degree of a bundle}
is $\sum k_j$.

Two collections $\{k_j x_j\}$ and
$\{l_j y_j\}$
are equivalent if their difference
$\{k_j x_j\} \cup\{-l_j y_j\}$
is the collection of zeros and poles of a meromorphic
function on  $\frC$.

In our case,
the line bundle  $\cL$ is of degree (-3).
Respectively, the dual bundle has degree
 $3$. Holomorphic sections of the dual bundle
can be obtained in the following way:
take a lineal functional on $\C^3$, it determines
a linear functional on each fiber of the bundle.
It turns to be that a section has 3 zeros
(taking in account multiplicities).
Next  (this is a result of calculations),
we set two points of the divisor to
two prescribed points, therefore the third point is
a parameter on the set of bundles.

After this we can write formulas for generators of
 $\GL(2,\Z)$ in these coordinates.

\smallskip

The author thanks A.~G.~Kuznetsov, A.~A.~Rosly, and S.~L.~Tregub
for discussion of this topic.

\section{Coordinates}

\COUNTERS

{\bf\punct Coordinates on the set of curves.}
We write the equation
$$
\det(\lambda+\nu A+\nu B)=0
$$
in the form
\begin{multline}
\lambda^3+ d_1 \mu^3+ d_2 \nu^3+
\\
+p_+ \lambda^2\mu+ p_-  \lambda \mu^2+\quad
q_+ \lambda^2\nu+ q_- \lambda \nu^2+\quad
r_+ \mu^2 \nu + r_-\mu\nu^2+
\\
+ t\lambda\mu\nu=0
\label{eq:krivaya}
\end{multline}

Reduce the matrix $A$ to diagonal form. Denote by
$h_1$, $h_2$, $h_3$ its eigenvalues.
 Note that  $-h_1$, $-h_2$, $-h_3$
are the coordinates $\lambda:1:0$ of intersection of the cubic curve
with the line  $\nu=0$. The numbers
$h_1$, $h_2$, $h_3$ are defined up to a permutation.
We fix their order.

  Denote by
$$
  U=\begin{pmatrix} u_{11}&u_{12}& u_{13}\\
                    u_{21}&u_{22}& u_{23}\\
                    u_{31}&u_{32}& u_{33}
  \end{pmatrix}
$$
the matrix of the operator $B$ in the eigenbasis of the matrix  $A$.
Notice that $U$ is defined up to a conjugation by a diagonal matrix.
For matrices in general position we can assume $u_{12}=u_{13}=1$.
Then all other coordinates
  $u_{ij}$
are fixed.

Therefore the total number of coordinates on the space of pairs
of matrices is
$3+7=10$.  The number of coefficients in the equation
(\ref{eq:krivaya}) is 9. Thus we need  at least one
additional coordinate
to reconstruct a pair of matrices from the spectral curve.

Obviously,
\begin{align}
d_1&=h_1h_2h_3,
\label{eq:d1}
\\
d_2&=\det(U)=\det(B),
\label{eq:d2}
\\
p_+&=h_1+h_2+h_3,
\label{eq:pplus}
\\
p_-&=h_1h_2+h_1h_3+h_2h_3,
\label{eq:pminus}
\\
q_+&=u_{11}+u_{22}+u_{33},
\label{eq:qplus}
\\
q_-&=\begin{vmatrix}u_{11}&u_{12}\\ u_{21}& u_{22}  \end{vmatrix}
   + \begin{vmatrix}u_{11}&u_{13}\\ u_{31}& u_{33}  \end{vmatrix}
   +  \begin{vmatrix}u_{22}&u_{23}\\ u_{32}& u_{33}  \end{vmatrix}
,
\label{eq:qminus}
\\
r_+&=h_1h_2 u_{33}+h_1h_3 u_{22}+h_2h_3 u_{11}
,
\label{eq:rplus}
\\
r_-&=h_3\begin{vmatrix}u_{11}&u_{12}\\ u_{21}& u_{22}  \end{vmatrix}
     + h_2\begin{vmatrix}u_{11}&u_{13}\\ u_{31}& u_{33}  \end{vmatrix}
     + h_1 \begin{vmatrix}u_{22}&u_{23}\\ u_{32}& u_{33}  \end{vmatrix}
,
\label{eq:rminus}
 \\
t&=(h_1+h_2)u_{33}+(h_1+h_3)u_{22}+(h_2+h_3)u_{11}
\label{eq:t}
.
\end{align}
Note that an (ordered) triple
 $d_1$, $p_+$, $p_-$
and an (unordered) triple
  $h_1$, $h_2$, $h_3$ determine one another.


\smallskip

{\bf\punct Partial reconstruction of the matrix
 $U$.}
Solving the system  (\ref{eq:qplus}),
(\ref{eq:rplus}), (\ref{eq:t}) of linear equations
for diagonal elements of the matrix   $U$ we get
\begin{align}
u_{11}=
\frac{q_+h_1^2-th_1+r_1}{(h_1-h_2)(h_1-h_3)}
,
\label{eq:u11}
\\
u_{22}=
\frac{q_+h_2^2-th_2+r_1}{(h_2-h_1)(h_2-h_3)}
,
\label{eq:u22}
\\
u_{33}=
\frac{q_+h_3^2-th_3+r_1}{(h_3-h_1)(h_2-h_2)}
.
\label{eq:u33}
\end{align}


{\bf\punct Divisor and coordinates of the bundle.}
The first coordinate on the space
$\C^3$ can be regarded as a section
of the bundle  $\cL^*$ on $\frC$.
Let us find zeros of this section.
In other words, we must find matrices
 $\lambda+\mu A+\nu B$ whose kernel contains
a vector with first coordinate $=0$, i.e.,
the following equation has a nonzero solution
  $v_2$, $v_3$:
$$
\begin{pmatrix}\lambda+\mu h_1+ \nu u_{11}&\nu  u_{12}&\nu  u_{13}\\
 \nu     u_{21}&\lambda+\mu h_2+ \nu u_{22}&\nu  u_{23}\\
  \nu     u_{31}&\nu  u_{32}&\lambda+\mu h_3+\nu  u_{33}
  \end{pmatrix}
\begin{pmatrix}
0\\ v_2\\ v_3
\end{pmatrix}
=
\begin{pmatrix}
 0\\0\\0
\end{pmatrix}
.
$$
Thus, we get    the following condition for
 $\lambda:\mu:\nu$
\begin{equation}
\rk
\begin{pmatrix}\nu  u_{12}&\nu  u_{13}\\
                \lambda+\mu h_2+ \nu u_{22}&\nu  u_{23}\\
                \nu  u_{32}&\lambda+\mu h_3+\nu  u_{33}
  \end{pmatrix}
=
1
\label{eq:rk23}
.
\end{equation}
There is an obvious pair of solutions
$$
\lambda:\mu:\nu = h_2:(-1):0,\qquad h_3:(-1):0
\,.$$
Equating minors of the matrix
 (\ref{eq:rk23}) to 0, we come to
\begin{align}
& \nu^2 u_{12}u_{23}-\nu u_{13}(\lambda+\mu h_2+ \nu u_{22})=0
,
\\
&\nu u_{12}(\lambda+\mu h_3+\nu u_{33})-\nu^2 u_{13} u_{32}=0
,
\\
& (\lambda+\mu h_2+\nu u_{22}) (\lambda+\mu h_3+\nu u_{33})
-\nu^2 u_{23} u_{32}=0
.
\end{align}
Assuming
 $\nu=1$ and solving two equation
with respect to $\lambda$, $\mu$,
we get the third point of the divisor,
$\lambda:\mu:\nu=L:M:1$, where
\begin{align}
  L&=\frac{u_{12}
 h_3\begin{vmatrix}u_{12}& u_{13}\\u_{22}&u_{23} \end{vmatrix}
+u_{13} h_2
\begin{vmatrix}u_{12}&u_{13}\\ u_{32}&u_{33}  \end{vmatrix}
 }
{u_{13}u_{12} (h_3-h_2)}
.
\\
  M&=-\frac{u_{12}
 \begin{vmatrix}u_{12}& u_{13}\\u_{22}&u_{23} \end{vmatrix}
+u_{13}  \begin{vmatrix}u_{12}&u_{13}\\ u_{32}&u_{33}  \end{vmatrix}
 }
{u_{13}u_{12} (h_3-h_2)}
.
\end{align}

Automatically, the point
 $L:M:1$ lies on the curve  $C$ (by the condition
(\ref{eq:rk23})).

\begin{theorem}
a) The point $L:M:1$
is uniquely determined by
a pair of matrices  $(A,B)\in\frM(3)$
and an ordering of eigenvalues of
a matrix  $A$.

\smallskip

b) A pair of matrices  $(A,B)\in\frM(3)$
of general position can be reconstructed
by the spectral curve and the point $L:M:1$.
\end{theorem}

The first statement is obvious, the collection
$$
(h_2:(-1):1),\qquad (h_3:(-1):1),\qquad (L:M:1)
$$
is a divisor determining a linear bundle of degree 3
on the elliptic curve. Since two points of divisor are fixed,
the third one is uniquely determined by a linear bundle and
determines a bundle.

The inverse construction is given by explicit formulas,
which are given in the next subsection.

\smallskip

{\bf\punct The inverse map.} Recall that
$u_{11}$, $u_{22}$, $u_{33}$
were reconstructed above  (\ref{eq:u11})--(\ref{eq:u33}).
Next, we can set
\begin{equation}
  u_{12}=1,\qquad u_{13}=1
.
\label{eq:u12}
\end{equation}
Then the expressions for $L$, $M$ come to a simpler form:
\begin{align*}
L&=\frac{h_3(u_{23}-u_{22})+h_2(u_{33}-u_{32})}{h_3-h_2}
,\\
M& =\frac{(u_{23}-u_{22})+(u_{33}-u_{32})}{h_3-h_2}
.
\end{align*}
We can regard this pair of equalities
 as a system of equations for
 $(u_{23}-u_{22})$ and $(u_{33}-u_{32})$.
Solving it, we get
\begin{align}
u_{23}=L+h_2 M+ u_{22} = L+h_2 M+
 \frac{q_+h_2^2-th_2+r_+}{(h_2-h_1)(h_2-h_3)}
,
\label{eq:u23}
\\
u_{32}=L+h_3 M+ u_{33} = L+h_3 M+
\frac{q_+h_3^2-th_3+r_+}{(h_3-h_1)(h_2-h_2)}
.
\label{eq:u32}
\end{align}
It remains two unknowns
 $u_{21}$, $u_{31}$ and two non-used equations
  (\ref{eq:qminus}), (\ref{eq:rminus});
  the equations are linear with respect
to unknowns. Theorem is proved (a pair of matrices can be
reconstructed from spectral data), but we need final expressions,
\begin{multline}
u_{21}=\frac{h_1-h_2}{h_2-h_3}(L+h_2M)(L+h_3M)+
\\
+
\frac{M}{(h_1-h_3)(h_2-h_3)}
\Bigl(r_+(h_1-h_2-h_3)-q_+ h_1h_2h_3+th_2h_3\Bigr)
+\\+
\frac{L}{(h_1-h_3)(h_2-h_3)}
\Bigl( q_+(h_2h_3-h_1h_2-h_1h_3)-r_++th_1   \Bigr)
+\\+
\frac{r_--q_-h_2}{h_2-h_3}+
\frac{(q_+h_2^2-th_2+r_+)(q_+h_1^2-th_1+r_1)}
{(h_1-h_2)^2(h_3-h_1)(h_2-h_3)}
\label{eq:u21}
,\end{multline}
\begin{multline}
u_{31}=\frac{h_1-h_3}{h_3-h_2}(L+h_2M)(L+h_3M)+
\\
+
\frac{M}{(h_1-h_2)(h_3-h_2)}
\Bigl(r_+(h_1-h_2-h_3)-q_+ h_1h_2h_3+th_2h_3\Bigr)
+\\+
\frac{L}{(h_1-h_2)(h_3-h_2)}
\Bigl( q_+(h_2h_3-h_1h_2-h_1h_3)-r_++th_1   \Bigr)
+\\+
\frac{r_--q_-h_3}{h_3-h_2}+
\frac{(q_+h_3^2-th_3+r_+)(q_+h_1^2-th_1+r_1)}
{(h_1-h_3)^2(h_2-h_1)(h_3-h_2)}
\label{eq:u31}
.\end{multline}
These formulas are
 cumbersome but not senseless. Observe, for instance,
elementary symmetric functions of variables  $(-h_1)$, $h_2$, $h_3$.
See also cumbersome expressions in the next Section,
 (\ref{eq:wtr1})--(\ref{eq:wtr4})
and (\ref{eq:long}).

\section{Action of generators of $\GL(2,\Z)$}

\COUNTERS

{\bf\punct  Transposition of matrices.}
Transposition $(A,B)\mapsto(B,A)$
corresponds to the following reflection of
the plane
$$
  \lambda:\mu:\nu \mapsto \lambda:\nu:\mu
.$$
The bundle on the curve is the same.

Our coordinates are not completely convenient for description
of this operation. However we can describe it in the terms of
the usual  'constructions with a cubic and a ruler'
(see, e.g., \cite{Cle}).

Denote by
 $\xi_1$, $\xi_2$, $\xi_3$ the eigenvalues
of the matrix
$B$. We use the following
notation for points of the curve
 $\frC$:
\begin{align*}
 P_1=(h_1:(-1):0),\quad   P_2:=(h_2:(-1):0),\quad
  P_3:=(h_3:(-1):0),\quad
\\
X_1:=(\xi_1:0:(-1)),\quad X_2:=(\xi_2:0:(-1)),\quad
X_3:=(\xi_3:0:(-1))
.
\end{align*}
Denote also
$$
Q=(L:M:1)
,$$
thus points  $\{P_2$, $P_3$, $Q\}$ are a divisor of
the bundle
$\cL^*$. We must find
the equivalent divisor of the form
 $X_2$, $X_3$, $Y$
with an unknown point  $Y$.

For this purpose consider the line
 $X_1 Q$, let $\ell_1(\lambda,\mu,\nu)=0$
be its equation. Let  $T$ be the third
point of intersection of this line and the cubic
 $\frC$.

Next, consider the line
 $P_1 T$, let   $\ell_2(\lambda,\mu,\nu)=0$ be its equation.

{\it The third  point of intersection
of the line  $P_1 T$ with the cubic $\frC$
is the  desired point $Y$.}

To prove this, we consider the following meromorphic
function on
 $\C\P^2$:
$$
f(\lambda:\mu:\nu)=\frac \mu\nu \cdot\frac{\ell_2}{\ell_1}
.
$$
Let us restrict this function to the curve
 $\frC$. It is easy to see that:

\smallskip

a) $P_2$, $P_3$, $Q$ are poles of the function $f$;

\smallskip

b) $X_2$, $X_3$ and the point $Y$ constructed
just now are zeros of $f$;

\smallskip

c) $X_1$, $P_1$, $T$ are removable singularities.

\smallskip

Therefore the divisors
 $P_2$, $P_3$, $Q$
and $X_2$, $X_3$, $Y$ are equivalent.

\begin{proposition}
The matrix  $\begin{pmatrix} 0&1\\1&0  \end{pmatrix}\in\GL(2,\Z)$
corresponds to the following transformation:

\smallskip

--- we transpose  $\mu$ and $\nu$ n the equation;
(\ref{eq:krivaya})

\smallskip

-- the distinguished point of this curve is the image of the point
$Y$ under the transposition of the coordinates $\mu$, $\nu$.
\end{proposition}


{\bf\punct The inversion of the first matrix.}
We consider the operation
 $(A,B)\mapsto (A^{-1},B)$.
The equation of new curve is given by
\begin{equation}
\det(\lambda+\mu A^{-1}+\nu B) =0
\label{eq:pas}
\end{equation}
or
$$
 \begin{vmatrix}
 \lambda+\mu h_1^{-1}+ \nu u_{11}&\nu  u_{12}&\nu  u_{13}\\
 \nu     u_{21}&\lambda+\mu h_2^{-1}+ \nu u_{22}&\nu  u_{23}\\
 \nu     u_{31}&\nu  u_{32}&\lambda+\mu h_3^{-1}+\nu  u_{33}
  \end{vmatrix}
.
$$
The new coordinates
 (see the notation (\ref{eq:krivaya}) are
\begin{align*}
\wt d_1&=h_1^{-1}h_2^{-1}h_3^{-1},
\\
\wt d_2&=\det(U)=\det(C),
\\
\wt p_+&=h_1^{-1}+h_2^{-1}+h_3^{-1},
\\
\wt p_-&=h_1^{-1}h_2^{-1}+h_1^{-1}h_3^{-1}+h_2^{-1}h_3^{-1}
,
\\
\wt q_+&=u_{11}+u_{22}+u_{33}
,
\\
\wt q_-&=
\begin{vmatrix}u_{11}&u_{12}\\ u_{21}& u_{22}  \end{vmatrix}
   + \begin{vmatrix}u_{11}&u_{13}\\ u_{31}& u_{33}  \end{vmatrix}
   +  \begin{vmatrix}u_{22}&u_{23}\\ u_{32}& u_{33}  \end{vmatrix}
,
\\
\wt r_+&=
h_1^{-1}h_2^{-1} u_{33}+h_1^{-1}h_3^{-1} u_{22}+h_2^{-1}h_3^{-1} u_{11}
,
\\
\wt r_-&=
h_3^{-1}\begin{vmatrix}u_{11}&u_{12}\\ u_{21}& u_{22}  \end{vmatrix}
+ h_2^{-1}\begin{vmatrix}u_{11}&u_{13}\\ u_{31}& u_{33}  \end{vmatrix}
+ h_1^{-1} \begin{vmatrix}u_{22}&u_{23}\\ u_{32}& u_{33}  \end{vmatrix}
,
 \\
\wt t&=
(h_1^{-1}+h_2^{-1})u_{33}+(h_1^{-1}+h_3^{-1})u_{22}+(h_2^{-1}+h_3^{-1})u_{11}
.
\end{align*}
Also, we must substitute
 $h_j\mapsto h_j^{-1}$
to the expression for $L$, $M$.

\begin{proposition}
The matrix $\begin{pmatrix} -1&0\\0&1 \end{pmatrix}\in\GL(2,\Z)$
corresponds to the following
transformation
of the space of cubic curves
with the distinguished points
\begin{align}
\wt d_1&=d_1^{-1}
,
\nonumber
\\
\wt d_2&= d_2
,
\nonumber
\\
\wt p_+&=\frac{p_-}{d_1}
,
\nonumber
\\
\wt p_-&=\frac{p_+}{d_1}
,
\nonumber
\\
\wt q_+&=q_+
,
\nonumber
\\
\wt q_-&=q_-
,
\nonumber
\\
\wt r_+&=\frac{q_+p_+ - t}{d_1}
,
\nonumber
\\
\wt r_-&=
-\frac1{h_1 h_2 h_3} \Bigl[
(h_1-h_2)(h_1-h_3)(L+M h_2)(L+M h_3)+
\label{eq:wtr1}
\\ &\qquad
+L\bigl(r_+- t h_1+ q_+(-h_2 h_3+h_1h_3+h_1 h_2)\bigr)+
\\ &\qquad
+M\bigr(r_+(h_2+h_3-h_1)-t h_2 h_3+q_+ h_1h_2h_3\bigr)+
\\ &\qquad
+q_- h_1(h_2+h_3)-h_1 r_-\Bigr]
\label{eq:wtr4}
\\
\wt t&=\frac{q_+ p_- - r_+ }{d_1}
,
\nonumber
\\
\wt L&=L+M(h_2+h_3)
,
\nonumber
\\
\wt M&=-h_2h_3 M.
\nonumber
\end{align}

\end{proposition}

All lines of this table are sufficiently obvious except $r_-$ which
requires a long calculation with
a usage of formulas of
inversion (\ref{eq:u11})--(\ref{eq:u33}) и
(\ref{eq:u12})--(\ref{eq:u32}).


\smallskip

{\bf\punct The third generator of $\GL(2,\Z)$.}
Nest, we consider the transformation
 $(A,B)\mapsto (A, AB)$,
and respectively the coefficients
of the equation
$$
  \det(\lambda+\mu A+\nu AB)=0
.
$$
This equation can be written in the form
$$
  \det(\lambda A^{-1}+\mu+\nu B)=0
.$$
This is equivalent to
 (\ref{eq:pas}) up to a change of notation.


\begin{proposition}
The matrix  $\begin{pmatrix}1&1\\0&1 \end{pmatrix}$
corresponds to the following
transformation of the space of curves with
distinguished points: the coefficients
$d_1$, $p_\pm$ are the same, other
coefficients transforms as
\begin{align}
d_2^*&= d_1 d_2,
\nonumber
\\
q_+^*&=q_+ p_+-t,
\nonumber
\\
q_-^*&=d_1 \wt r_-,
\label{eq:long}
\\
r_+^*&=d_1 q_+,
\nonumber
\\
r_-^*&=d_1 q_-,
\nonumber
\\
t^*&= p_-q_+-r_+,
\nonumber
\\
L^*&=-h_2h_3 M,
\nonumber
\\
M^*&=L+M(h_2+h_3),
\nonumber
\end{align}
where $\wt r_-$ is given by  (\ref{eq:wtr1})--(\ref{eq:wtr4})
\end{proposition}

{\tt Math.Dept., University of Vienna,

 Nordbergstrasse, 15,
Vienna, Austria

\&


Institute for Theoretical and Experimental Physics,

Bolshaya Cheremushkinskaya, 25, Moscow 117259,
Russia

\&

Mech.Math. Dept., Moscow State University,
Vorob'evy Gory, Moscow


e-mail: neretin(at) mccme.ru

URL:www.mat.univie.ac.at/$\sim$neretin

wwwth.itep.ru/$\sim$neretin
}

\end{document}